\documentclass[a4paper]{article}

\usepackage{latexsym}
\usepackage{amsfonts}

\newtheorem{definition}{Definition}
\newtheorem{theorem}{Theorem}

\newcommand{\setR}{\ensuremath{\mathbb{R}}}
\newcommand{\rot}{\ensuremath{\mathrm{rot}}}
\newcommand{\dD}{\mathrm{d}}

\begin{document}

\title{Star Products for integrable Poisson Structures on $\setR^3$}
\author{      Claus Nowak\\
Fakult\"at f\"ur Physik der Universit\"at Freiburg\\
Hermann-Herder-Str. 3, 79104 Freiburg i.Br. / Germany}
\date{August 1997\\[2mm] 
      Freiburg Preprint THEP 97/14}

\maketitle

\begin{abstract}
We prove the existence of a deformation quantization for 
integrable Poisson structures on $\setR^3$ and give a generalization 
for a special class of three dimensional manifolds.
\end{abstract}

\vspace{2cm}

The program of deformation quantization of the function algebra
on a symplectic manifold extends naturally to manifolds with 
nonregular Poisson structures.
In contrast to symplectic manifolds the existence of star 
products on nonregular Poisson manifolds, even on $\setR^n$, is 
an open problem. Particular examples of quantizable nonregular 
Poisson structures were found, e.g. a star product for linear 
Poisson structures \cite{Dri83,Gutt83} which is induced by a star 
product on the cotangent bundle of a Lie group, 
and for quadratic Poisson structures in three 
dimensions \cite{DoninMakar95}. We will give in this 
letter a proof for the existence of star products for integrable Poisson 
structures on $\setR^3$ and extend this to a class of three dimensional
manifolds.

On $\setR^3$ the components $P^{ij}$ of a Poisson structure 
$P=\frac12 P^{ij}\partial_i\wedge\partial_j$ 
can be arranged as a three-dimensional Poisson vector 
$\vec{P}=(P^{23},P^{31},P^{12})$. The Jacobi identity for $P$ is 
equivalent to the identity $\vec{P}(\rot\vec{P})=0$. Solutions of 
this equation are of the form $\vec{P}=\psi\vec{\nabla}\varphi$ 
for arbitrary functions $\psi$ and $\varphi$. In the following 
we will concentrate on 
\begin{equation}\label{defPvektor}
  \vec{P}=\vec{\nabla}\varphi
\end{equation}
which are called {\em integrable} Poisson structures \cite{AlbeFei}. 

The main result of this letter is stated in the following theorem:
\begin{theorem}
  Let $P$ be an integrable Poisson structure on $\setR^3$. Then 
  there exists a star product which is a quantization of this 
  Poisson structure.
\end{theorem}

{\bf Proof:}
We will construct a star 
product of Weyl type. We are looking for product which is a formal 
power series of bilinear operators $M_k$, $k\geq 0$, vanishing on 
constants for $k\geq 1$,
such that
\begin{displaymath} 
  f\ast g~:=~\sum_{k=0}^\infty \hbar^k M_k(f,g)
\end{displaymath}
is an associative product on $C^\infty(\setR^3)[[\hbar]]$ with 
$M_0(f,g)=fg$ and $M_1(f,g)=\frac{i\hbar}{2}\{f,g\}$. 
The product is of Weyl type if the $M_k$ 
have the parity 
\begin{equation}\label{parity}
  M_k(f,g)=(-1)^k M_k(g,f)~~.   
\end{equation}

\noindent Associativity is equivalent to \cite{LecomtedeWilde}
\begin{equation}\label{defHochbed}
  \delta(M_k)~=~\frac12 \sum_{l=1}^{k-1} [M_l,M_{k-l}]=:R_k
\end{equation}
with the Hochschild-$\delta$ for a multilinear map $M$ with
$m+1$ arguments 
\begin{eqnarray}\label{defdelta}
  \delta(M)(f_0,\ldots,f_{m+1})&:=&f_0M(f_1,\ldots,f_{m+1})\\\nonumber
 &&\mbox{\hspace{-3cm}} -\sum_{i=0}^{m}(-1)^iM(f_0,
      \ldots,f_{i}f_{i+1},\ldots,f_{m+1}) 
   +(-1)^m M(f_0,\ldots,f_{m})f_{m+1} \nonumber
\end{eqnarray}
the Gerstenhaber bracket ($N$ multilinear with $n+1$ arguments)
\begin{equation}\label{defGerstenhaberKlammer}
  [M,N]~:=M\circ N -(-1)^{mn} N\circ M
\end{equation}
and the Gerstenhaber product 
\begin{displaymath} 
  (M\circ N)(f_0,\ldots,f_{m+n})=\sum_{i=0}^m(-1)^{in}
         M(f_0,\ldots,N(f_i,\ldots,f_{i+n}),\ldots,f_{m+n})
\end{displaymath}
Since $\delta^2=0$ solving (\ref{defHochbed}) is a cohomological problem
\cite{LecomtedeWilde}. The necessary condition $\delta R_k$ is 
satisfied for $M_2\ldots M_{k-1}$ fulfilling (\ref{defHochbed}) 
due to the fact that the Gerstenhaber bracket is a graded Lie 
bracket. Since the cohomology of $\delta$ is nontrivial, we have to 
show that $R_k$ lies in the trivial cohomology class. 
The application of the theorem of Hochschild-Kostant-Rosenberg 
on the algebra of local multi differential operators on the function 
algebra of a manifold $M$ \cite{CahenGuttDeWilde}
gives the cohomology classes of $\delta$:
\begin{displaymath} 
  H^k_{\mathrm{Hoch}}(M)~\cong~\Gamma(\Lambda^{k+1}TM)
\end{displaymath}
The different indices $k$ for the cohomology classes and $k+1$ for 
the multi vector fields is a consequence of the grading of the 
operators which is the number of arguments minus one.
We will need this theorem for local multilinear maps which are 
differential operators in each argument, i.e. which vanish on constant
functions, the cohomology is the same, however, for the complex 
of all local multilinear maps \cite{CahenGuttDeWilde}.

The decomposition $R_k=\delta(M_k)+\alpha_k$ with $\alpha_k\in
\Gamma(\Lambda^{3}TM)$ together with the fact that $\delta(M_k)$
is a sum of maps which are symmetric in two consecutive 
arguments leads to the identity 
\begin{displaymath}
  AR_k~=~\alpha_k
\end{displaymath}
where $AR_k$ is the antisymmetrization of $R_k$. This can be read in the 
following way: If $R_k$ is a multilinear map with $\delta(R_k)=0$ 
then the antisymmetrisation $A R_k$ is a totally antisymmetric 
differential operator which is 1-differential in each argument.

If $M_2\ldots M_{k-1}$ were constructed with the symmetry 
(\ref{parity}) then $R_k(f,g,h)=-(-1)^k R_k(h,g,f)$, 
i.e. $R_k$ is symmetric for $k$ odd and the condition $AR_k=0$
is trivially fulfilled. On the other hand for
$M(f,g)=\pm M(g,f)$ is $\delta(M)(f,g,h)=\mp\delta(M)(h,g,f)$ and 
therefore $M_k$ can be chosen with the right symmetry if the 
integrability condition $AR_k=0$ is fulfilled.

We will now focus on an important property of $\delta$. Writing  
out explicitly
\begin{displaymath} 
  M(f_0,\ldots,f_m)=\sum_{|I_0|,\ldots,|I_m|}M^{I_0,\ldots,I_m}
      \partial_{I_0}f_0\ldots\partial_{I_m}f_m
\end{displaymath}
with multiindices $I_0,\ldots,I_m$, 
$I_j=(i_{j,1},\ldots,i_{j,|I_j|})$,
$\partial_{I_j}=\partial_{i_{j,1}}\ldots 
\partial_{i_{j,|I_j|}}$, and Einstein summation convention for all 
indices of the multiindices,
then 
\begin{eqnarray*}
  (\delta M)(f_0,\ldots,f_{m+1})&=& \\
 && \mbox{\hspace{-3.2cm}} -\sum_{i=0}^{m}(-1)^i 
  \sum_{|I_0|,\ldots,|I_m|}\sum_{J\cup K=I_i}
  M^{I_0,\ldots,I_m}~\cdot  \partial_{I_0}f_0\ldots\partial_{J}f_i
   \partial_{K}f_{i+1}\ldots\partial_{I_m}f_{m+1}  
\end{eqnarray*}
where $J\cup K=I_i$ is the sum over all partitions of the multiindex 
$I_i$. We see: The coefficients of the multidifferential 
operator $\delta(M)$ are linear combinations of the coefficients 
of $M$. We call this a combinatorial operator. By examination of the 
proof of the Hochschild-Kostant-Rosenberg theorem in 
\cite{CahenGuttDeWilde} one sees the inverse statement:
If $N=\delta(M)$, then there exists an $\tilde{M}$ with 
$N=\delta(\tilde{M})$ such that the
coefficients of $\tilde{M}$ are linear combinations of the coefficients 
of $N$. 

We will recursively prove the existence of bilinear differential 
operators $M_k$ which satisfy (\ref{defHochbed}) and have the 
following additional properties:
\begin{enumerate}
  \item $M_k$ is an ordered polynomial operator (This is defined below). 

  \item The total number of derivatives of the $M_k$ is greater 
  than two for $k\geq 2$.
\end{enumerate}
\begin{definition}
  A multidifferential operator $M$ is called an ordered polynomial 
  operator (OPO), if the coefficients of $M$ are polynomials in  
  the coefficients of Poisson structure and in partial derivatives 
  of the coefficients which can be arranged such that the indices of 
  the Poisson 
  structures are contracted with partial derivatives which stand 
  on the right side of the Poisson structure. 
\end{definition}
Examples:
\begin{itemize}
  \item The Poisson bracket $\{f,g\}=P^{ij}\partial_if\partial_jg$ 
  is an OPO

  \item The terms of the Jacobi identity $\{\{f,g\},h\}$ are OPO.

  \item The operator $\partial_rP^{is}\partial_sP^{jr}\partial_if 
  \partial_jg$ is \underline{not} an OPO.
\end{itemize}

It is clear that the property of being ordered is well defined for 
polynomial coefficients since for each term there are only a finite 
number of arrangements of the factors, if one of these arrangements 
has the required property then the term is an OPO, if none of the 
arrangements has the property then it is not an OPO.

We emphasize that this property is a property of the operator with 
the given coefficients, it is unfortunately not a property 
of the operator itself. This is due to the fact that there 
exist operators which can be written as OPO or not as OPO. 
An example is given by the operator 
\begin{equation}\label{exOPOnonOPO}
  \partial_kP^{ij}\partial_i(P^{kr}\partial_rP^{lm}+
  P^{lr}\partial_rP^{mk}+P^{mr}\partial_rP^{kl})\partial_j\partial_l 
  f \partial_m g~=~ 0
\end{equation}
which vanishes due to the Jacobi identity. 
There is one term in the expansion of 
(\ref{exOPOnonOPO}), $\partial_kP^{ij}P^{kr}\partial_i\partial_rP^{lm}
\partial_j\partial_l f \partial_m g$, which is an OPO, all other 
terms are not OPO. Using (\ref{exOPOnonOPO}) we can express this 
OPO as an operator which is not an OPO.

We will now list a few properties of the Gerstenhaber bracket
(\ref{defGerstenhaberKlammer}) and the Hochschild-$\delta$
(\ref{defdelta}):
\begin{itemize}
  \item  $M,~N$ OPO $\Rightarrow [M,N]$ OPO. This is clear since 
  in the Gerstenhaber product one operator is inserted in the 
  other which preserves the structure.

  \item $M$ OPO $\Rightarrow \delta(M)$ OPO. This is a consequence 
  of the combinatorial property of $\delta$.

  \item $N=\delta(M), N$ OPO $\Rightarrow$ it is possible to choose 
  $M$ as an OPO. This is a consequence of the combinatorial 
  property of the inverse of $\delta$.
\end{itemize}
For $M_1,\ldots,M_{k-1}$ OPO is $R_k$ an OPO. 
If it is possible to show that $AR_k=0$, then it is clear that 
$M_k$ can be chosen as an OPO. Since $R_k$ is a differential operator
with 3 arguments it has the total differential degree 
$\geq 3$, therefore $M_k$ will be a differential operator of degree 
$\geq 3$. 

We are now ready to prove $AR_k=0$ for $k$ even. $AR_k$ is a 
sum of totally antisymmetric differential operators of degree $(1,1,1)$
(where the degree of an operator is the sequence of the 
differential degrees in each argument)
These terms emerge from terms in 
$R_k=\frac12 \sum_{l=1}^{k-1}[M_l,M_{k-l}]$ 
which are of degree $(1,1,1)$ (since the antisymmetrization does not 
change the degree) This  
is only possible, if one of $M_l$ or $M_{k-l}$ is a differential 
operator of degree $(1,1)$, this must be $M_1$, since all other $M_l$ 
have a total degree $\geq 3$. Therefore 
\begin{eqnarray} \nonumber
  AR_k(f,g,h) &=& (A(M_{k-1}\circ M_1)(f,g,h))_{\mbox{1-diff}}\\
  && \mbox{\hspace{-2cm}}=~\frac13(M_{k-1}(M_1(f,g),h)-M_{k-1}(f,M_1(g,h))
              ~+~\mbox{cycl.})_{\mbox{1-diff}} \nonumber \\
    &&\mbox{\hspace{-2cm}}=~\frac23(M_{k-1}(M_1(f,g),h)  
              ~+~\mbox{cycl.})_{\mbox{1-diff}}  \label{reducedCycl}
\end{eqnarray}
We have to take explicitly the 1-differential parts (indicated by 
the subscript 1-diff) since leaving 
out terms in $R_k$ destroys the cocycle property and therefore the 
property of the antisymmetrization of being 1-differential. 
The terms of degree $(1,1,1)$ in $M_{k-1}(M_1(f,g),h)$ come from the 
terms $M_{k-1}^{(m,1)}$ which are of degree $(m,1)$. These are of 
the form
\begin{displaymath}
  M_{k-1}^{(m,1)i_1\ldots i_m j}\partial_{i_1}\ldots\partial_{i_m} f
  \partial_j g
\end{displaymath}
$M_{k-1}$ is an OPO, therefore the $M_{k-1}^{(m,1)i_1\ldots i_m j}$ are 
polynomials in coefficients of the Poisson structure and partial 
derivatives thereof. The indices of the last Poisson 
structure of this ordered series have to be contracted with 
partial derivatives standing more right, i.e. partial derivatives 
of the arguments. Since the indices of the Poisson structure are 
antisymmetric, it must be partial derivatives of different 
arguments. Therefore $j$ must be one of the indices, the other is one of  
$i_1,\ldots,i_m$. Without loss of generality we will assume 
that $i_m$ is this index and write now
\begin{displaymath}
  M_{k-1}^{(m,1)}(f,g)~=~M_{k-1}^{(m,1)I i j}\partial_I\partial_{i} f
  \partial_j g
\end{displaymath}
where $I=(i_1,\ldots,i_{m-1})$ is a multiindex, $M_{k-1}^{(m,1)I i 
j}$ is antisymmetric in $i$ and $j$.
The part of degree $(1,1,1)$ in (\ref{reducedCycl}) coming from 
$M_{k-1}^{(m,1)}$ is then, up to a factor, given by
\begin{equation}\label{restCycl}
   M_{k-1}^{(m,1)I i j}\partial_I\partial_{i} P^{rs}\partial_r f
   \partial_s g\partial_j h +\mbox{cycl.} 
\end{equation}
We will show that these terms vanish for each $m$ separately, 
for this we will now use the fact that the space is $\setR^3$ and 
the special form (\ref{defPvektor}) of the Poisson structure.

Since (\ref{restCycl}) is derivative in its arguments 
it is sufficient to use coordinate functions to show the vanishing of 
(\ref{restCycl}). It is trivially fulfilled if two of the three 
functions $f,g$ and $h$ are equal. It remains to check the 
condition for $f=x^1,~g=x^2$ and $h=x^3$:

\begin{eqnarray*}
&& \mbox{\hspace{-1cm}}M_{k-1}^{(m,1)I i 3}\partial_I\partial_{i} P^{12}+
M_{k-1}^{(m,1)I i 1}\partial_I\partial_{i} P^{23}+
M_{k-1}^{(m,1)I i 2}\partial_I\partial_{i} P^{31} \\
  &=&M_{k-1}^{(m,1)I 1 3}\partial_I\partial_{1}\partial_{3}\phi+
     M_{k-1}^{(m,1)I 2 3}\partial_I\partial_{2}\partial_{3}\phi+
     M_{k-1}^{(m,1)I 2 1}\partial_I\partial_{2}\partial_{1}\phi\\
  &&  + M_{k-1}^{(m,1)I 3 1}\partial_I\partial_{3}\partial_{1}\phi+
     M_{k-1}^{(m,1)I 1 2}\partial_I\partial_{1}\partial_{2}\phi+ 
     M_{k-1}^{(m,1)I 3 2}\partial_I\partial_{3}\partial_{2}\phi \\
  &=& 0
\end{eqnarray*}
because of the symmetry of the partial derivatives and the 
antisymmetry of $M_{k-1}^{(m,1)I i j}$ in $i$ and $j$, from line 1 to 
line 2 we have written out the summation over $i$ explicitly. 

Hence we conclude $AR_k=0$ and we can continue the recursion with an
operator with the desired properties. This ends the proof of the 
existence of a star product for the considered Poisson structures.

It ist clear that it is difficult to transfer this proof on other 
Poisson structures. For this special form of the Poisson structure 
the Jacobi identity $\vec{P}(\rot\vec{P})=0$ is fulfilled because 
$\rot\vec{P}=0$, this seems to be the reason for $AR_k=0$. 
The proof of the existence of a star product on $\setR^3$ for the general 
case $\vec{P}=\psi\vec{\nabla} \varphi$ in this manner fails, 
an explicit calculation shows that it is necessary to add a non OPO to 
$M_3$ in order to get $AR_4=0$.

There exist a generalization for special type of three dimensional 
manifolds given in the following theorem:
\begin{theorem} 
Let $M$ be an orientable, three dimensional 
manifold with a flat torsion free connection and a covariant
constant volume form $\mu=\frac16\mu_{ijk}\dD x^i\wedge\dD x^j\wedge\dD 
x^k$. Define an isomorphism
\begin{eqnarray}\nonumber
  \ast: \Gamma(\Lambda^2 TM) &\to & \Gamma(\Lambda^1 T^\ast M) \\
        Y &\mapsto & i(Y)\mu \nonumber
\end{eqnarray}
with the injection $i(Y)\mu=\frac12Y^{ij}\mu_{ijk}\dD x^k$ for a bivector
$Y=\frac12 Y^{ij}\partial_i\wedge\partial_j$.
Then for all $\alpha\in\Gamma(\Lambda^1 T^\ast M)$ with $\dD\alpha=0$ the 
bivector $\ast^{-1}\alpha$ is a Poisson bivector and there exists a star 
product for this Poisson structure.
\end{theorem}

{\bf Proof:} For a flat torsion free connection it is always
possible to find a chart around each point such that the Christoffel 
symbols vanish.  
In such a chart the coefficients $\mu_{ijk}$ of the covariant constant 
volume form are constant and proportional to $\epsilon_{ijk}$, without 
loss of 
generality we will assume that $\mu_{ijk}=\epsilon_{ijk}$. We will use in 
the following an atlas containing only such charts.
For a one-form
$\alpha=\alpha_i\dD x^i$ the bivector $\ast^{-1}\alpha$
is given by $P^{ij}=\epsilon^{ijk}\alpha_k$ 
and the Jacobi identity for $P$ is then 
\begin{eqnarray*}
  P^{1r}\partial_rP^{23}\!\!\!\!&+&\!\!\!\!P^{2r}\partial_rP^{31}
        +P^{3r}\partial_rP^{12} 
  ~ = ~ P^{1r}\partial_r\alpha_1+P^{2r}\partial_r\alpha_2
        +P^{3r}\partial_r\alpha_3 \\
  & = &  P^{12}(\partial_2\alpha_1-\partial_1\alpha_2)+ 
     P^{23}(\partial_3\alpha_2-\partial_2\alpha_3)+      
     P^{31}(\partial_1\alpha_3-\partial_3\alpha_1) \\
  & = & 0
\end{eqnarray*}
for $\dD\alpha=0$, i.e. $\ast^{-1}\alpha$ is a Poisson bivector.

In the domain of a chart the closed one-form can be written as
$\alpha=\dD\phi$, i.e. the Poisson structure is locally integrable.
Using the arguments for $\setR^3$ we 
conclude that there exists a star product in this chart, i.e. the 
condition $AR_k=0$ is fulfilled. But for $M_1,\ldots,M_{k-1}$ 
given on the whole manifold, forming $R_k$ by (\ref{defHochbed}) is 
a coordinate invariant construction, and the vanishing of the total 
antisymmetrisation is also a coordinate invariant property. 
Therefore $AR_k=0$ ist fulfilled in every chart and hence on the whole 
manifold and $M_k$ can be constructed globally, it has the required 
properties in every chart. 
\vspace{5mm}

\noindent{\bf Acknowledgement}

\noindent I would like to thank Martin Bordemann for pointing out 
the generalization given in the last theorem.

\end{document}